\documentstyle[fleqn]{article}
\font\BBb=msbm10 at 12pt
\newcommand{\Bbb}[1]{\mbox{\BBb #1}}

\newcommand{\be}{\begin{equation}}
      \newcommand{\ee}{\end{equation}}
      \newcommand{\ba}{\begin{eqnarray}}
       \newcommand{\ea}{\end{eqnarray}}
\newcommand{\ban}{\begin{eqnarray*}}
\newcommand{\ean}{\end{eqnarray*}}

 \newcommand{\qed}{\hspace*{\fill}\rule{3mm}{3mm}\quad \vspace{.2cm}}
 \newcommand{\Pf}{\noindent {\bf Proof:} }

\newcommand{\sect}[1]{\section{#1} \setcounter{equation}{0}}

\newtheorem{theo}{Theorem}[section]
\newtheorem{defn}{Definition}[section]

\begin{document}
\newtheorem{lem}[theo]{Lemma}
\newtheorem{prop}[theo]{Proposition}  
\newtheorem{coro}[theo]{Corollary}
\newtheorem{ex}[theo]{Example}
\newtheorem{note}[theo]{Note}

\title{Hausdorff Convergence and Universal Covers }
\author{Christina Sormani
  \and Guofang Wei\thanks {Partially supported by NSF Grant \# DMS-9971833.}}
\date{}
\maketitle

\begin{abstract}
We prove that if $Y$ is the Gromov-Hausdorff limit of 
a sequence of compact manifolds, $M^n_i$, with
a uniform lower bound on Ricci curvature and a uniform upper bound
on diameter, then $Y$ has a universal cover.
We then show that, for $i$ sufficiently large, the fundamental
group of $M_i$ has a surjective homeomorphism  onto the group of deck 
transforms of $Y$.  Finally, in the non-collapsed
case where the $M_i$ have an additional uniform lower bound on volume, we 
prove that the kernels of these surjective maps are finite with a uniform
bound on their cardinality.  A number of theorems are also proven 
concerning the limits of covering spaces and their deck transforms
when the $M_i$ are only assumed to be compact length spaces with a 
uniform upper bound on diameter.
\end{abstract}

\newcommand{\inj}{\mbox{inj}}
\newcommand{\vol}{\mbox{vol}}
\newcommand{\diam}{\mbox{diam}}
\newcommand{\Ric}{\mbox{Ric}}
\newcommand{\Iso}{\mbox{Iso}}
\newcommand{\Hess}{\mbox{Hess}} 
\newcommand{\divg}{\mbox{div}}
\newcommand{\RR}{\bf{R}}

\sect{Introduction}

In recent years the
limit spaces of manifolds with lower bounds on Ricci curvature
have been studied from both a geometric and topological perspective.
In particular, Cheeger and Colding have proven a number of results regarding
the regularity and geometric properties of these spaces.  However, the
topology of the limit spaces is less well understood.  Note that in this paper
a manifold is a Riemannian manifold without boundary.

Anderson \cite{An} has proven that there are only finitely 
many isomorphism types
of fundamental groups of manifolds with a uniform upper bound on diameter,
lower bound on volume  and lower bound on
Ricci curvature.  Thus one might think that given a converging sequence 
of such manifolds, the fundamental groups of the manifolds must
eventually be isomorphic to the fundamental group of the limit space.
However, Otsu \cite{Ot} has shown that there are metrics of uniformly
positive Ricci curvature on $S^3 \times {\Bbb R}P^2$ 
which converge to a simply connected 5-dim metric space, showing 
that this need not be the case.  

Tuschmann \cite{Tu} has proven that if $Y$ is 
the limit space of a sequence of manifolds
with two sided sectional curvature bounds then $Y$ is locally simply 
connected and thus has a universal cover.  
In fact Perelman \cite{Pl} shows that 
the limit space of a sequence of manifolds
with a lower bound on sectional curvature is locally contractible. 
If the limit 
space is locally simply connected, it is not difficult to show that 
eventually there is a surjective map from the
fundamental groups of the manifolds onto the fundamental group of the limit
space (see \cite{Tu, Ca}, \cite[Page 100]{Gr}, also Section 2 of this paper).
Zhu  \cite{Zh} has proven a similar result for limits of three dimensional
manifolds with uniform lower bounds on Ricci curvature and volume and a uniform
upper bound on diameter. 

In this paper a manifold is a Riemannian manifold without boundary.
Here we are concerned with limits of sequences of manifolds with
a uniform upper bound on diameter and lower bound on Ricci curvature.
The limits of such sequences have only been shown to be locally 
simply connected at special ``regular''  points \cite{ChCo}.  In fact
Menguy \cite{Me} has 
shown that the limit space could locally have infinite topological type.

We prove that the universal cover of  the limit space exists.  
We can thus study the group of deck transforms of the universal cover, 
$\bar{\pi}_1(M)$ [Defn~\ref{deck}].  
Note that this revised fundamental group, $\bar{\pi}_1(M)$, 
is isomorphic to the fundamental
group of $M$ if $M$ is locally simply connected (c.f. [Sp]).  
We can now state the main theorem of our paper.

\begin{theo} \label{surjective} 
Let $M_i$ be a sequence of compact manifolds satisfying
\be \label{cond1}
Ricci(M_i) \ge (n-1)H \textrm{ and } diam(M_i) \leq D,  
\ee
for some $H \in {\Bbb R}$ and $D >0$.
If $Y$ is the Gromov-Hausdorff limit of the $M_i$ 
then the universal cover of $Y$ exists and for $N$ sufficiently large
depending on $Y$, 
there is a surjective homeomorphism
\be
\Phi_i:\pi_1(M_i) \to \bar{\pi}_1(Y) \qquad \forall i \ge N.
\ee
\end{theo}

\begin{note}
We don't known if $Y$ is semi-locally simply connected, or the universal 
cover is simply connected yet. See Section 2, 
Example~\ref{sin}, \ref{ring} for more information.
\end{note}

When the sequence is non-collapsing we have a stronger result:

\begin{theo} \label{noncol}
Let $M_i^n$ be a sequence of compact manifolds satisfying 
\be \label{cond2}
Ricci(M_i) \ge (n-1)H, \qquad diam(M_i) \leq D \ \textrm{ and } vol(M_i)\ge V  
\ee
for some $H \in {\Bbb R}$, $D>0$ and $V>0$.
If $Y$ is the Gromov-Hausdorff limit of the $M_i$
then there is $i_0 = i_0(n,v,D,\delta_Y)$ such that 
$\pi_1(M_i)/F_i$ is isomorphic to $\bar{\pi}_1 (Y)$ for all $i \geq i_0$,
 here $F_i$ is a finite subgroup of $\pi_1(M_i)$, and the order of each 
 $F_i$ is uniformly bounded by $N(n,v,D)$. In particular, 
 $\pi_1(M_i)/F_i$ is isomorphic to
  $\pi_1(M_j)/F_j$ for all $i,j \geq i_0$,
\end{theo}

Compare Anderson's result \cite{An} which says that there 
are only finitely many isomorphism types
of fundamental groups of compact manifolds satisfying (\ref{cond2}).

To prove these results we need to study the limits spaces of 
compact length spaces.  Thus in Sections 2 and 3 we restrict
ourselves to  sequences, $M_i$, which are only {\em compact length spaces} 
with $diam(M_i) \le D$ that converge in the
Gromov-Hausdorff sense to a limit space $Y$. 

In Section 2 we present two examples of such sequences of length spaces 
which converge in the Gromov-Hausdorff topology.
However, their fundamental groups
cannot be mapped surjectively onto the fundamental group or revised
fundamental group of the limit space.         
In the first example, we have a sequence of simply connected
length spaces whose limit space is not simply connected [Ex~\ref{sin}].
In the second example,
the limit space has no universal cover [Ex~\ref{ring}]. 
It should be recalled that even when the limit space
is a manifold that the limit of the universal covers is not
necessarily a cover of the limit space (see \cite[Theorem 2.1]{Pe1} 
for a case where it is).  Thus the universal
cover cannot be directly used to prove properties about the
fundamental group.

In Section 3, we introduce $\delta$-covering spaces [Defn~\ref{defdel}].
Unlike the universal cover, $\delta$-covers always exist.    
We then show that 
the limit of the $\delta$-covers is a cover of the limit space 
[Theorem~\ref{limdelcov}].  
Furthermore, we prove that for a fixed $\delta>0$
groups of deck transforms of the $\delta$-covers, $\tilde{M}^\delta_i$,
of the $M_i$ eventually have a surjective map onto
deck transforms of the $\delta$-cover of the limit space [Cor~\ref{GG}].
We also describe
the relationship with the $\delta$-covers and the universal cover if 
latter exists [Theorem~\ref{exists}].
We conclude with a proof of the following theorem which should
be contrasted with Example~\ref{sin} [Note~\ref{nsimpconn}].

\begin{theo} \label{simpconn}
If $M_i$ are simply connected compact length spaces 
with $diam(M_i) \le D$ that converge
in the Gromov Hausdorff topology to a compact space $Y$, 
then the universal cover of $Y$ exists and the revised
fundamental group, $\bar{\pi}_1(Y)$,
is trivial.
\end{theo}
 
In Section 4 we study 
limit spaces, $Y$, of sequences of 
compact manifolds $M_i$ satisfying (\ref{cond1})
and prove that the universal cover of $Y$ exists.
To do so we first prove Lemma~\ref{LemSnk} 
extending techniques from [So] involving the Excess Theorem 
of [AbGl] to arbitrary lower bounds on Ricci curvature.  
We then apply this lemma combined with
[ChCo] regularity results to  prove Theorem~\ref{regdel} 
that there is at least one ball in
the limit space that lifts isometrically to any covering space.  
To extend this to balls around arbitrary points in the limit
space we use the Bishop-Gromov Volume Comparison, which holds
for limit measures by [ChCo] and Theorem~\ref{limdelcov}. 
Thus we prove that there is a $\delta_Y>0$ such
that the $\delta_Y$-cover of the limit space $Y$ is the universal cover of $Y$
[Theorem~\ref{samedelta}].
Note that this $\delta_Y$ depends on many factors and cannot be determined 
uniformly without contradicting the examples of Otsu [Ot].
Combining Theorem~\ref{samedelta} with the results of Section 3, 
we obtain Theorem~\ref{surjective} and Corollary~\ref{theoisom}. 

We conclude with the non-collapsed case, where $M_i$ satisfy
(\ref{cond2}), proving
Theorem~\ref{noncol}. We use results from Section 3 
regarding the kernel of the surjective map and results of Anderson 
regarding elements of the fundamental group of manifolds with lower 
bounds on volume [An] to obtain Theorem~\ref{noncol}.

Background material for Gromov-Hausdorff limits and
Ricci curvature can be found in 
Chapter 1 Sections A-C, Chapter 3 Sections A-B,  
Chapter 5 Section A of \cite{Gr} and in Chapters 9-10 of \cite{Pe2}.  
Background material
on covering spaces and fundamental groups can be found in
Chapters 1-2 of \cite{Sp} and  \cite{Ma}.

%---------------------------------------------------------------------------

\sect{Background and Examples}

%----------------------------------------------------------------------------

In Sections 2 and 3 we consider compact length spaces.  
No curvature condition is assumed.
See \cite[Chapter 1]{Gr} for basic results about length spaces (called 
path metric spaces).  Recall also (c.f. \cite[Chapter 3A]{Gr}), the
following definition of the Gromov-Hausdorff distance between metric spaces
(called Hausdorff distance).

\begin{defn} {\em \cite[Defn 3.4]{Gr}
Given two metric spaces $X$ and $Y$, the {\em Gromov-Hausdorff distance}
between them is defined,
\be
d_{GH}(X,Y)=\inf \left\{ d^Z_H(f(X),f(Y)): \stackrel
{ \textrm{ for all metric spaces }Z, \textrm{ and isometric } }
{ \textrm{embeddings }f: X \rightarrow Z, g: Y \rightarrow Z }
\right\},
\ee
where, $d_H^Z$ is the Hausdorff distance between subsets of $Z$,
\be
d_H^Z(A,B)=\inf\{\epsilon>0: B \subset T_\epsilon(A) \textrm{ and }
A \subset T_\epsilon(B) \}.
\ee
Here $T_\epsilon(A)=\{x\in Z: d_Z(x,A) <\epsilon\}$.

If $x\in X$ and $y\in Y$, the {\em pointed Gromov-Hausdorff distance}:
\be
d_{GH}((X,x), (Y,y))=\inf \left\{ d^Z_H(f(X),f(Y)): \stackrel
{ \textrm{ for all metric spaces } Z, \textrm{ and isometric embed-} }
{ \textrm{dings }f: X \rightarrow Z, g: Y \rightarrow Z 
\textrm{ s.t.} f(x)=g(y) }
\right\}.                               
\ee }
\end{defn}

It is then clear what we mean by the Gromov-Hausdorff convergence of 
compact metric spaces.  However, for non-compact metric spaces, 
the following looser definition of convergence is used.

\begin{defn} {\em \cite[Defn 3.14]{Gr} 
We say that non-compact length spaces $(X_n, x_n)$ converge in 
the Gromov-Hausdorff sense to $(Y,y)$ if 
for any $R>0$ there exists a sequence $\epsilon_n \to 0$ such that 
$B_{x_n}(R+\epsilon_n)$ converges to $B_{y}(R)$ in the Gromov-Hausdorff
sense.}
\end{defn}

Although we are limiting ourselves to compact length spaces, their
universal covers may well be non-compact.

Recall that a space, $Y$, is {\em semi-locally simply connected} (or
semi-locally one connected) if for all $y\in Y$ there is a neighborhood
$U$ of $y$ such that $\pi_1(U, y) \to \pi_1(Y,y)$ is trivial 
(\cite[p 78]{Sp} \cite[p 142]{Ma}).  
That is, any curve in $U$ is contractible in $Y$.
This is weaker than saying that $U$ is simply connected.

For a metric space $Y$, let $r(Y)$ denote the maximal number $r$ such that 
every 
closed curve in a ball of radius $r$ in $Y$ is homotopic to zero 
in $Y$, the semi-locally simply connectivity radius. Note 
that if $Y$ is compact and semi-locally simply connected, then $r(Y)$ is 
positive.   

The following theorem demonstrates how Gromov Hausdorff closeness affects the
fundamental groups of spaces which are semi-locally simply connected.

\begin{theo}  \label{theosurj}
Let $Y_1, Y_2$ be two compact length spaces and $Y_2$ is semi-locally 
simply connected. If $d_{GH} (Y_1, Y_2) \leq \epsilon$, where $\epsilon = 
r(Y_2)/20$, then there exists a surjective homomorphism $\Phi:\ \pi_1(Y_1) 
\rightarrow  \pi_1(Y_2)$.
\end{theo}

 This theorem essentially follows from the proof of \cite[Theorem(b)]{Tu}. 
See also \cite[Page 100]{Gr} and \cite{Ca}. As an extension of this theorem, 
Theorem~\ref{theodelsurj}, will be proven in the next section,
 we will omit the proof.

\begin{note} 
One usually doesn't get a similar surjection for high homotopy groups, 
as Berger's  examples of $S^3$ collapsing to $S^2$ illustrate for
$\pi_2$. 
\end{note}

Theorem 2.1 immediately gives the following isomorphism.
 
\begin{coro}
Let $Y_1$ and $Y_2$ be two compact semi-locally simply connected length spaces
and let  $r = \min \{ r(Y_1)/20,r(Y_2)/20 \}$.
 If $d_{GH} (Y_1, Y_2) \leq r$ then $\pi_1(Y_1)$ is isomorphic to 
 $\pi_1(Y_2)$.
\end{coro}

\begin{note}{\em  Note that if $Y_1=S^2\times S^1_{\epsilon}$ and
$Y_2=S^2$ then $\pi_1(Y_1)$ only maps surjectively onto
$\pi_1(Y_2)$.  This does not contradict the above corollary
because $r$ there depends on both $r(Y_1)=\epsilon$ and $r(Y_2)$.}
\end{note}

\begin{note} Note that Colding \cite{Co} and Cheeger-Colding  \cite{ChCo}, 
proved that if 
$Y_1,Y_2$ are same dimensional manifolds and sufficiently Hausdorff close 
(closeness depends on $Y_1, Y_2$), then $Y_1,Y_2$ are diffeomorphic. 
Here an assumption on $r(Y_i)$ is shown to suffice to match
the fundamental groups of length spaces.
\end{note}

Recall the definition of the universal cover of a metric space 
[Sp, page 62,80]. First, $Y$ is a covering space of $X$ if there is
a continuous map $\pi: Y \to X$ such that $\forall x \in X$
there is an open neighborhood $U$ such that $\pi^{-1}(U)$ is a disjoint union
of open subsets of $Y$ each of which is mapped homeomorphically onto $U$
by $\pi$.   We say that connected space $\tilde{X}$ is a universal 
cover of $X$ if $\tilde{X}$
is a cover of $X$ such that for any other cover $Y$ of $X$, there is a 
a commutative triangle formed by a continuous map $f:\tilde{X} \to Y$
and the two covering projections.  

Recall that if $\pi: E \rightarrow Y$ is a covering and $Y$ is a connected 
compact length 
space, then there is a unique length metric on $E$ making 
$\pi: E \rightarrow Y$  
distance non-increasing and a local isometry (see \cite{Ri}).
Of course, the covering space need not be compact.

The universal cover may not exist as can be seen in 
[Sp, Ex 17, p 84].  However, if it exists, then it is unique.  Furthermore,
if a space is locally path connected and semi-locally simply connected 
then it has a universal cover and that cover is simply connected
[Sp, Cor 14, p 83].  On the other hand, the universal covering space
of a locally path connected space may not be simply connected
[Sp, Ex 18, p 84].  We now present two examples.

\begin{ex} \label{sin}
There exists a compact length space, $X$, which is the
limit of simply connected compact length spaces, but
is not simply connected.
\end{ex}

\Pf
The construction of $X$ is based upon the fact that
$\{(t, \sin (1/t)): 0<t\leq 1/\pi\} \cup \{0\}\times [-1,1]$
is a compact metric space which isn't path connected.
Here, however, we construct a path connected space
which isn't simply connected.

First define two compact
 sets:
 $$
  K_1=\{(x, y, \sin(1/y)): x\in [-1/\pi,1/\pi], |y| \le |x| \}
 $$
 $$
  K_2=\{(x,y,z):  |x| =|y| \le 1/\pi,  z\in [-2,2]\}
 $$

 Let the length space be
 $$
 X= K_1 \cup K_2.
 $$                                           
First we prove that $X$ is not simply connected.

We claim that the loops $C_j:[0,4/(j\pi)]\to X$ defined as follows are
 all homotopic to each other but are not contractible.
 $$
 C_j(t)=(t,-t,0) \subset K_2 \qquad t\in [0, 1/(j\pi)]
 $$
 $$
 C_j(t)=(1/(j\pi), t - 2/(j\pi), 0)\subset K_1
\qquad t\in [1/(j\pi), 3/(j\pi)]
 $$        
 $$
C_j(t)=(4/(j\pi) - t, 4/(j\pi) - t, 0) \qquad t \in [3/(j\pi), 4/(j\pi)]
 $$
 Clearly the length of this curve is
 $$
 L(C_j)=(2\sqrt{2}+2)/(j\pi) \to 0, \qquad j\to \infty
 $$
 and clearly each curve is homotopic to the next.

 Now why aren't they contractible?
 If $C_1$ were contractable, there would be a homotopy   
 $$
 H: A=[0,1/\pi]\times[-1/\pi,1/\pi] \to X \subset R^3
 $$
 such that
 $$
 H(0,t)=(0,0,0), H(s,-1/\pi)=(s, -s, 0), H(s,1/\pi)=(s,s,0),
 H(1/\pi,t)=(1/\pi,t,0).
 $$

Write $H(s,t)=(h_1(s,t), h_2(s,t), h_3(s,t))$.

It is easy to see that $X_r=X - K_1 \cap  \{(x,y,z): z > r\}$   
is not simply connected for $r=3/4$, and that $C_1$ is not contractible
in this set.  Thus $H^{-1}(X_{3/4})$ is also not simply
connected and thus cannot contain the entire domain.
So there exists $(s_1, t_1)\in A$ such that
$h_3(s_1, t_1) > 3/4$ and $H(s_1,t_1) \in K_1$.

Let $U_1 \subset A$
be the connected component of $H^{-1}(X \setminus X_{1/2})$ that contains
$(s_1, t_1)$.
                  
We continue this definition inductively, given
$(s_j, t_j)$ we define $U_j \subset A$ to be the
connected component of $H^{-1}(X \setminus X_{1/2})$ that contains
$(s_j, t_j)$.  Then we note that $C_1$ is not contractable in
$H^{-1}(X_{3/4}\cup H(U_1) \cup ... \cup H(U_j))$ because at most
finitely many peaks have been filled in.  So there exists
$(s_{j+1}, t_{j+1})\in A \setminus ( U_1 \cup ...\cup U_j)$
such that $h_3(s_{j+1}, t_{j+1})> 3/4$ and $H(s_{j+1}, t_{j+1}) \in K_1$.

Now $A$ is compact so a subsequence of the $(s_j, t_j)$ converges to   
some $(s_\infty, t_\infty)\in A$.  By continuity
$h_3(s_\infty, t_\infty)\ge 3/4$.  Now  each point in the
subsequence is in a distinct connected component of
$H^{-1}(X \setminus X_{1/2})$.  Thus there exist points in $H^{-1}(X_{1/2})$
on line segments between the points in the subsequence.
The limit of these points must have $h_3 \le 1/2$.
 Since these line segments are shorter and
shorter the limit point of points on the segments must also
be $(s_\infty, t_\infty)$ and we have a contradiction.
No homotopy can exist.     

Now we define $X_k$, a sequence of simply connected length
spaces which converge to $X$.
$$
K_{k,1}=\{(x, y, \sin(1/y)):
     1/(k\pi) \le |x| \le 1/\pi, 1/(k\pi) \le |y| \le |x| \}
$$
$$
K_{k,1}'=\{(x, y, 0):  |x| \le 1/(k\pi), |y| \le |x| \}
$$
and the connecting set:
$$
K_{k,2}=\{(x,y,z):  |x| =|y| \le 1/\pi,  z\in [-2,2]\}
$$

Then the length spaces are defined,
$$
  X_k= K_{k,1} \cup K'_{k,1} \cup K_{k,2}.
$$            

It is easy to see that these spaces are compact with
a uniform bound on their diameters.  They are
also easily seen to be simply connected.

Now we claim $d_{GH}(X_k, X) < 10/(k\pi)$.
Let $Z_k$ be a metric space defined as
the disjoint union of $X_k$ and $X$ with
$K_2$ identified with $K_{k,2}$
and $K_{k,1}$ identified with the appropriate
part of $K_1$.  It is easy to see that
the tubular neighborhood of radius $10/(k\pi)$
about $X_k$ in $Z_k$ contains $X$ and visa versa
because the tubular neighborhood about $K_2$
covers $Z_k \setminus (X \cup X_k)$.
\qed

\begin{ex} \label{ring}
There exists a compact length space, $X$, with no universal cover
which is the limit of a sequence of compact length spaces, $X_i$,
that have universal covers.  Furthermore, the fundamental
groups of the $X_i$ are finitely generated but the fundamental group
of the limit space is not.
\end{ex}

\Pf
Let $X$ be the Hawaii Ring.  That is $X\subset R^2$ is defined
$$
X=\bigcup_{j=1}^\infty C_j
$$
where $C_j$ is the circle of radius $1/j$ around $(1/j, 0)$.

Let $X_i=\bigcup_{j=1}^i C_j$.
\qed

If a space has a universal cover, then we can define a revised
fundamental group.  First recall the definition of the group of
deck transforms (or self equivalences) of a cover $\pi: Y \to X$.
This is the group of homeomorphisms $h: Y \to Y$ such that
$\pi \circ h = \pi$ [Sp p85].  It is denoted $G(Y, X)$.

Given $\tilde{p} \in Y$, there is a natural surjection $\Psi_{\tilde{p}}$
from the fundamental group, $\pi_1(X, p)$ to $G(Y, X)$ defined as follows.
Note that given $g\in \pi_1(M,p)$ and $\tilde{p}$, 
then it has a representative loop, $C$.  We can lift $C$ to a curve based
at $\tilde{p}$ in the cover.  This defines an action of $g$ on 
$\pi^{-1}(p)\subset \tilde{M}$ which can be extended uniquely to a deck 
transform of $M$.  This map is surjective when $Y$ is path connected
because given
any $h \in G(Y, X)$ we can join $\tilde{p}$ to $h(\tilde{p})$ by a curve
and then the curve's projection to the base space gives an element 
of $\pi_1(X, p)$.  The kernel, $H_p$,
 consists of elements of the fundamental group $\pi_1(M, p)$,
whose representative loops are still
closed when they are lifted to the cover.

Note that when the universal cover
is simply connected and locally path connected, then the group
of deck transforms of the universal cover is isomorphic to the fundamental
group of $X$ [Sp page 87 Cor 4]. 

\begin{defn}  \label{deck} {\em
Let the {\em revised fundamental group} of $M$, denoted
 $\bar{\pi}_1(M)$, be the group of deck transforms of the universal
cover of $M$.}
\end{defn}

Then $\bar{\pi}_1(M)=\pi_1(M, p)/H_p$ where $H_p$ is a normal subgroup
described above.  It
is isomorphic to $\pi_1(M, p)$ when $M$ is semi locally simply connected.
In fact $H_p$ consists of those elements of $\pi_1(M, p)$ whose representative
loops lift to closed curves in the universal cover.  Thus if 
$\bar{\pi}_1(M, p)=\pi_1(M, p)$ then $H_p$ is trivial, so every
loop in $M$ which lifts to a closed curve in $\tilde{M}$ is contractible.
In particular, each loop in $\tilde{M}$ must project to a closed loop
in $M$ which is contractible, so it must be contractible.  Thus
$\bar{\pi}_1(M)=\pi_1(M, p)$ iff the universal cover of $M$ is 
simply connected.

\begin{note} \label{nsimpconn}
In Example~\ref{sin}, we have a sequence of simply connected spaces which
converge to a space which is not simply connected.  However, the limit
space is its own universal cover.  Thus its {\em revised fundamental group}
is trivial.  In Section 3 we will prove that this is always the case
[Theorem~\ref{simpconn}]. 
\end{note}

%---------------------------------------------------------------------------
\sect{$\delta$-Covering Spaces}
%-------------------------------------------------------------------------

In this section we continue to study locally path connected length spaces.
We would like to understand what happens to the fundamental groups of
such spaces when we take their limits.  Since we do not assume that they
are semi-locally simply connected, we do not know if these spaces have
simply connected universal covers.  In fact they might not have universal 
covers at all.  To circumvent this problem we introduce intermediate
$\delta$-covering spaces.

Let $\cal U$ be any open covering of $Y$. For any $p \in Y$, by 
\cite[Page 81]{Sp}, there is a covering space, $\tilde{Y}_{\cal U}$, 
of $Y$ with covering group  $\pi_1(Y,{\cal U}, p)$, where 
$\pi_1(Y,{\cal U}, p)$ is a normal subgroup of $\pi_1(Y, p)$, generated 
by homotopy classes of closed paths having a representative of the form 
$\alpha^{-1} \circ \beta \circ \alpha$, where $\beta$ is a closed path 
lying in some element of $\cal U$ and $\alpha$ is a path from $p$ to 
$\beta(0)$.  

\begin{defn} \label{defdel} {\em
Given $\delta>0$, the {\em $\delta$-cover}, denoted $\tilde{Y}^\delta$,
of a length space $Y$, is defined to be $\tilde{Y}_{{\cal U}_{\delta}}$
where ${\cal U}_\delta$ is the open covering of $Y$ consisting of
all balls of radius $\delta$.

The covering group will be denoted $\pi_1(Y,\delta, p)\subset \pi_1(Y,p)$
and the group of deck transforms of $\tilde{Y}^\delta$ will be denoted
$G(Y,\delta)=\pi_1(Y, p)/\pi_1(Y,\delta, p)$. }
\end{defn}

Note that $\tilde{Y}^{\delta_1}$ covers $\tilde{Y}^{\delta_2}$ when 
$\delta_1 \leq \delta_2$. 

Note also that $G(Y, \delta)$, denoted $G(\tilde{Y}^\delta, Y)$ in
[Sp, p 86 Cor 3], does not depend on $p$.
One can think of $G(Y, \delta)$ as roughly corresponding to the long loops,
of length at least $\delta$, in $\pi_1(Y,p)$.

There is a natural surjective map from $\pi_1(Y,p)$ to $G(Y, \delta)$, which
depends on $\tilde{p}\in \tilde{Y}^\delta$, 
that we can call $\Psi_{\tilde{p},\delta}$.  
(See the paragraph after Example~\ref{ring} in Section 2).

\begin{defn} \label{defnlength}
When $Y$ is compact, then
for any $h \in G(Y,\delta)$, we can define the (translative) 
$\delta$-length of $h$,
\be   
l(h, \delta)= \min_{q\in \tilde{Y}^\delta} d_{\tilde{Y}^\delta}(q, h(q)).
\ee
For any $g \in \pi_1(Y, p)$
\be   
l(g, \delta)= 
\min_{q\in \tilde{Y}^\delta} 
d_{\tilde{Y}^\delta}(q, \Psi_{\tilde{p},\delta}(g)(q)).
\ee
\end{defn}

Note that since this is a compact length space, there is a point 
$q\in \tilde{Y}$ and a curve $C$ from $q$ to $ \Psi_\delta(g)(q))$
whose length is $l(g, \delta)$.  Consequently there is a point $\pi(q)\in Y$
and a loop $\pi(C)$ based at the point, $q$, of length $l(g, \delta)$.
Furthermore there is some path, $\alpha$, from $p$ to $\pi(q)$ such that
$\alpha^{-1} (\pi(C)) \alpha$ is a representative curve for $g\in \pi_1(Y,p)$.

We have the following basic properties for $\delta$-length.

\begin{lem} \label{lemdel}
For all nontrivial $g \in G(Y, \delta,p)$, the $\delta$-length of $g$, 
\be
l(g, \delta) \geq \delta.  \label{lemdeleq}
\ee
For all $\delta_1 \leq \delta_2$ we have 
\be
l(g,\delta_1) \geq l(g,\delta_2). \label{lemdeleq2}
\ee
\end{lem}

\Pf
Suppose there exists a $g$ such that (\ref{lemdeleq}) fails to hold.
Since $Y$ is a length space, there is a curve $\beta$ of length
less than $\delta$ running from some point,
$q\in \tilde{Y}^\delta$ to $g(q)$.  There is also a curve
$\alpha$ running from $\tilde{p}$ to $q$.  Then $\pi(\alpha)$ runs 
$p$ to $\pi(q)$, and $\pi(\beta)$ is a loop in $B_{\pi(q)}(\delta)$,
so
$$
g=
[(\pi(\alpha)^{-1}\circ  \pi(\beta)\circ \pi(\alpha)]\in G(Y,\delta, p)
$$
is trivial by the definition.

As mentioned above, for all $\delta_1 \leq \delta_2$, there exists 
$\alpha$ and $\beta$ in $\tilde{Y}^{\delta_1}$ such that
$L(\beta)=l(g,\delta_1)$ and 
$$
g=[(\pi(\alpha)^{-1}\circ  \pi(\beta)\circ \pi(\alpha)]\in \pi(Y, p).
$$
Now $\tilde{Y}^{\delta_1}$ is the $\delta_1$ cover of $\tilde{Y}^{\delta_2}$,
so $\alpha$ and $\beta$ can be projected to paths in 
$\tilde{Y}^{\delta_1}$.  Thus $l(g,\delta_1) \geq l(g,\delta_2)$.
\qed

\begin{prop}  \label{uniisdel}
If $Y$ is a compact length space that has a universal cover,
then the universal cover is a $\delta$-cover, for some $\delta_0>0$.
Thus $\tilde{Y}=\tilde{Y}^{\delta_0}=\tilde{Y}^\delta$ for
all $\delta<\delta_0$.
\end{prop}

\Pf
If $Y$ has a universal cover, then for all $y\in Y$, there is
an $r_y>0$ such that $B_y(r_y)$ is lifted isometrically to
the universal cover. Suppose the universal cover is not a $\delta$-cover
for any $\delta>0$.  Then there exists $p_i\in Y$ and $C_i$
a loop based at $p_i$ in $B_{p_i}(1/i)$ which lift non-trivially  to the 
universal cover.  $Y$ is compact
so a subsequence of $p_i$ converge to $p_\infty$, and eventually
some $C_i$ will be contained in $B_{p_\infty}(r_{p_\infty})$ which
is a contradiction.
	\qed

\begin{prop}
If $Y$ is a compact length space such that the universal cover exists and all
nontrivial elements of $\pi_1(Y, p)$ have a positive $\delta$-length for some 
$\delta$,
then Y is semi-locally simply connected. 
\end{prop}

\Pf
Let $\delta_0$ be defined as in Proposition~\ref{uniisdel}.  If $Y$ is not
semi-locally simply connected, then there is a curve $C$ contained
in some $B_q(\delta_0) \subset Y$ which is not contractable in $Y$.
Let $\beta$ run from $p$ to $q$ and $g = [\beta \circ C \circ \beta^{-1}]
\in \pi_1(Y, p)$.  Note that $\Psi_{\delta_0}(g)$ is trivial
and since $\tilde{Y}^\delta=\tilde{Y}^{\delta_0}$ for all
$\delta<\delta_0$ by our choice of $\delta_0$, $\Psi_{\delta_0}(g)$
is trivial for all $\delta<\delta_0$.  This contradicts the hypothesis.
\qed

Note that in the example in Section 2, there is a nontrivial element
of the fundamental group with $0\ \delta$-length for all $\delta$.

Now we have good covering spaces for $Y$, but we don't know if 
Bishop-Gromov volume comparison holds on these covers even if it 
does on $Y$. In order to find a good covering space such that 
Bishop-Gromov volume comparison also holds, we relate these 
covering to the covering from the sequence. 
Note that in general,
the universal covering of sequence doesn't convergence to the universal 
cover of the limit, e.g. the lens spaces, $S^3/Z_p$, converges to $S^2$.  
(See \cite[Theorem 2.1]{Pe1} for a case 
where it does). 
However, we will show that this is almost true in the $\delta$-cover level 
(see Theorem~\ref{limdelcov}). To prove this first we need a 
revised version of Theorem~\ref{theosurj} which does
not require a bound on the semi-locally simply connectivity radius.  

\begin{theo}  \label{theodelsurj}
Let $Y_1, Y_2$ be two compact length spaces and such that 
$d_{GH} (Y_1, Y_2) \leq \epsilon$, 
then there is a surjective homomorphism,
$\Phi: G(Y_1,\delta_1) \to G(Y_2,\delta_2)$ 
for any $\delta_1>  20 \epsilon$ and $\delta_2>\delta_1+10\epsilon$.
\end{theo}

In the proof of this theorem and subsequent theorems,
we think of $G(Y, \delta)=\pi_1(Y, p)/\pi_1(Y,\delta, p)$ 
as $G(Y, \delta, p)$ consisting of equivalence classes of loops
based at $p$ such that $C_1 \sim C_2$ iff $C_1 \circ C_2^{-1}$
is homotopic to a curve in $\pi_1(Y, \delta, p)$.

\Pf 
Since $d_{GH} (Y_1, Y_2) \leq \epsilon$, there must exist a
metric space $(Z, d)$, such that $Y_1$ and $Y_2$ are isometrically
embedded in $(Z,d)$ and the Hausdorff distance between them is
less than $2 \epsilon$ and $d_Z(p_1, p_2)<2\epsilon$.
If $\alpha \in G(Y_1, \delta_1, p_1)$ then
it can be represented by some rectifiable closed curve 
$\gamma : [0,1] \rightarrow Y_1$. On $\Gamma := \gamma ([0,1])$ now 
choose points $x_1, \cdots, x_m, x_i = \gamma (t_i)$, with $x_m = x_0=p_1$ and 
$0 =t_0 \leq t_1 \leq \cdot \leq t_m = 1$ such that for each 
$\gamma_i := \gamma |_{[t_i,t_{i+1}]}, i= 0, \cdots, m-1$, one has 
$d(x_i, x_{i+1}) \leq L(\gamma_i) < 5\epsilon$. We will refer to this as 
a $5\epsilon$-partition of $\gamma$. 
Since $d_{H}^Z (Y_1, Y_2) \leq 2\epsilon$, 
for each $x_i, i=0, \cdot, m-1$, we can choose points $y_i \in Y_2$ 
satisfying $d_Z(x_i,y_i) < 2\epsilon$ and set 
$y_m := y_0=p_2$. Connected $y_i$ 
to $y_{i+1}$ by a minimizing geodesic, $\bar{\gamma}_i, i=0, \cdots, m-1$
of length less than $9\epsilon$. 
This yields a closed curve $\bar{\gamma}$ in $Y_2$, consisting of 
$m$ minimizing segments and having $p_2$ as its base point. Now define 
\[
\Phi (\alpha) = \Phi ([\gamma]) := [\bar{\gamma}] \in G (Y_2, \delta_2, p_2).
\]

 First we need to verify that $\Phi$ is a well defined. 
Using the fact that $9\epsilon<\delta_1/2 < \delta_2/2$, 
one easily see that $[\bar{\gamma}]$ doesn't depend on the choice of 
minimizing geodesics $\bar{\gamma}_i$, nor  on the choice of 
points $y_i \in Y_2$, nor on the special partition $\{x_1, \cdots, x_m \}$ 
of $\gamma ([0,1])$. 
Moreover using additionally the uniform continuity of a 
homotopy one can see that $[\bar{\gamma}]$ only depends on the homotopy 
class of $\gamma$.  
It thus also easy to check that $\Phi$ is a homomorphism from
$\pi_1(Y_1, p_1)$ to $G(Y_2, \delta_2, p_2)$.
However $\alpha \in G(Y_1, \delta_1, p_1)$ not $\pi_1(Y_1, p_1)$.

Suppose $\gamma_1$ and $\gamma_2$ are both representatives of $\alpha
\in G(Y_1, \delta_1, p_1)$.  Then $\gamma_1 *\gamma_2^{-1}$ is
homotopic to a loop $\gamma_3$ generated by loops of the
form  $(\alpha * \beta) * \alpha^{-1}$, where $\beta$ is a closed path lying 
in a ball of radius $\delta_1$ and $\alpha$ is a path from 
$p_1$ to $\beta(0)$.  
So $[\bar{\gamma_1}]=[\bar{\gamma_3}] * [\bar{\gamma_2}]$ and
we need only show that $[\bar{\gamma_3}]$ is trivial in 
$G(Y_2, \delta_2, p_2)$.

In fact $\bar{\gamma_3}$ can be chosen  as follows.
The $y_i$'s corresponding to the $x_i$'s from the $\beta$ segments of 
$\gamma_3$  are all within $\delta_1+2 \epsilon$ of a common point
and the minimal geodesics between them are within 
$\delta_1+(2 +9/2) \epsilon < \delta_2$.  Furthermore, the  $y_i$'s 
corresponding to the $x_i$'s from the $\alpha$ and $\alpha^{-1}$ segments 
of the curve can be chosen to correspond.  Thus  $\bar{\gamma_3}$
is generated by loops of the
form  $(\alpha * \beta) * \alpha^{-1}$, where $\beta$ is a closed path lying 
in a ball of radius $\delta_2$ and $\alpha$ is a path from $p_2$ to $\beta(0)$.
So it is trivial.

Last, we need to show that $\Phi$ is onto. If 
$\bar{\alpha} \in G(Y_2, \delta_2, p_2)$, it can be represented by some 
rectifiable closed curve $\sigma$.  Choose now an $\epsilon$-partition 
$\{y_0, \cdots, y_m\}$ of $\sigma$ such that for all $i=0, \cdots, m-1$, 
one has  $L(\sigma|_{[t_i,t_{i+1}]}) < \epsilon, (y_i = \sigma(t_i))$, 
further corresponding points $x_i \in Y_1$ satisfying 
$d(x_i,y_i) < 2\epsilon$ and we can connect those points by minimizing 
curves in $Y_1$. This yields a piecewise length minimizing 
$\gamma: [0,1] \rightarrow Y_1$ with base point $x_0$, and because of 
$d(x_i, x_{i+1}) < 5\epsilon$ the curve $\gamma$ allows a 
$5\epsilon$-partition and $[{\gamma}]\in G(Y_1, \delta_1, p_1)$.
Now $\Phi([\gamma])=\bar{\alpha} \in G(Y_2, \delta_2, p_2)$
because $\Phi([\gamma])$ was shown above not to depend upon the
choice of the $y_i$ as long as $d_Z(x_i,y_i)<2 \epsilon$ and 
the $x_i$ were a $5\epsilon$-partition of a representative curve
$\gamma$.

Therefore $\Phi$ is surjective. 
\qed

Let $(M_i, p_i)$ be a sequence of connected locally path-connected spaces 
that converge to $(Y,p)$ in the pointed Gromov-Hausdorff topology.  Denote 
$\tilde{M}_i^\delta$ the $\delta$-covering of $(M_i,p_i)$, 
$\pi_1(M_i,\delta, p_i)$ the covering group and 
$G(M_i,\delta, p_i) = \pi_1 (M_i, p_i)/\pi_1 (M_i, \delta, p_i)$ the deck 
transformation on $\tilde{M}_i^\delta$. 

Now Theorem~\ref{theodelsurj} gives us the following.

\begin{coro}  \label{GG}
If $(M_i, p_i)$ is a sequence of connected locally path-connected spaces with
diam$(M_i) \leq D$ which converges to $(Y,p)$ in the pointed Gromov-Hausdorff 
topology, 
then for any $\delta_1<\delta_2$, there exists $N$ sufficiently large depending
upon $\delta_2$ and $\delta_1$ such that $\forall i\ge N$ there is 
a surjective map $\Phi_i:G(M_i,\delta_1)  \to G(Y,\delta_2)$.
\end{coro}

Using this we can show that in Gromov-Hausdorff limit the $\delta$-covering 
of the sequence converges to a cover which is almost the  $\delta$-covering 
of the limit. The fact that it isn't quite the  $\delta$-covering is seen
by the following simple example.  Take a sequence of flat tori of side lengths
$1$ by $(n-1)/(2n)$.  The $\delta=1/2$ cover of these tori are cylinders
since all loops of length $< 1/2$ are not unraveled.  However, they converge 
to a torus of side lengths $1$ by $1/2$ whose $\delta=1/2$ cover is 
Euclidean space, which is a cover of the limit cylinder.  However, for 
any $\delta_2<\delta$, the $\delta_2$-cover of the limit is a cylinder.
 
\begin{theo}  \label{limdelcov}
 If $M_i$ with diam$(M_i) \leq D$ converges to $Y$ in the 
Gromov-Hausdorff metric and the 
$\delta$-covering of $M_i$, $(\tilde{M}_i^\delta, \tilde{p}_i)$, converges in 
the pointed Gromov-Hausdorff metric to $({Y}^\delta, \tilde{p}_\infty)$, then 
$({Y}^\delta, \tilde{p}_\infty)$ is a covering space of $Y$, which is covered
by the $\delta$-cover of $Y$, $\tilde{Y}^\delta$.  Furthermore,
for all $\delta_2>\delta$, $Y^\delta$ covers $\tilde{Y}^{\delta_2}$.
So we have, covering projections mapping
$$
\tilde{Y}^\delta \to Y^\delta \to \tilde{Y}^{\delta_2} \to Y.
$$
\end{theo}

\Pf Let $\pi_i^\delta: \tilde{M}_i^\delta \rightarrow 
M_i$ be the covering map. It's distance decreasing by construction. After 
possibly passing to a subsequence it follows from a generalized version of the 
Arzela-Ascoli theorem (see e.g. \cite[Page 279, Lemma 1.8]{Pe2}) 
that $\pi_i^\delta: \tilde{M}_i^\delta \rightarrow M_i$ will 
converge to a distance decreasing map $\pi^\delta: Y^\delta \rightarrow Y$. 

First, since $\tilde{M}_i^\delta$ is the $\delta$-covering space, the covering 
map $\pi_i^\delta: \tilde{M}_i^\delta \rightarrow M_i$  must be an 
isometry on any ball 
of radius $< \delta$. As $\pi_i^\delta$ converges to $\pi^\delta$ 
this property must be 
carried over to $\pi^\delta: Y^\delta \rightarrow Y$. 
Hence $\pi^\delta: Y^\delta \rightarrow Y$ is a covering map. 

So we have three covering spaces of $Y$,
$\tilde{Y}^\delta$, $Y^\delta$ and $\tilde{Y}^{\delta_2}$.  
By the Unique Lifting Theorem \cite[Lemma 3.1, Page 123]{Ma} 
if $\tilde{Y}_1$ and $\tilde{Y}_2$ are covers of $Y$, then $\tilde{Y}_1$
covers $\tilde{Y}_2$ if every closed curve in $Y$ which lifts to a closed
curve in $\tilde{Y}_1$ also lifts to a closed curve in $\tilde{Y}_2$.  

Now if $C$ is a closed curve in $Y$ whose lift to $\tilde{Y}^{\delta}$ is  
closed, then it is homotopic to a curve
consisting of paths, loops within $\delta$-balls and reverse paths.  So
its lift to $Y^\delta$ is also closed since $\pi^\delta$ is an 
isometry on $\delta$-balls. Therefore $\tilde{Y}^\delta$ covers $Y^\delta$.

If $\delta_2>\delta$, we want to show
$Y^\delta$ covers $\tilde{Y}^{\delta_2}$.  Suppose not.  Then there is
a closed curve $C$ in $Y$ whose lift to $Y^\delta$ is closed but
whose lift to $\tilde{Y}^{\delta_2}$ is not a closed loop.

Since the lift of $C$ in $\tilde{Y}^{\delta_2}$ is not closed, 
$\Phi_{\delta_2}([C]) \in G(Y,\delta_2)$ is nontrivial. 
Using Corollary~\ref{GG},
we can find $N$ sufficiently large so that 
$\Phi_i:G(M_i,\delta) \to G(Y,\delta_2)$
is surjective.  In particular we can find curves $C_i$ which converge
to $C$ in the Gromov-Hausdorff sense, such that $\Phi_i([C_i])=[C]$.
Since, $[C_i]$ are nontrivial their lifts to $\tilde{M}_i^\delta$
run between points $\tilde{C}_i(0)\neq\tilde{C}_i(1)$.  Furthermore, by
Lemma~\ref{lemdel}, 
$$
d_{\tilde{Y}_i^\delta}(\tilde{C}_i(0), \tilde{C}_i(1))\ge \delta. 
$$
In the limit, the lifted curves $\tilde{C_i}$ converge to the
lift of the limit of the curves, $\tilde{C}$ in $Y^\delta$  and
$$
d_{Y^\delta}(\tilde{C}(0), \tilde{C}(1))\ge \delta.  
$$
This implies that $\tilde{C}$ is not closed and we have a contradiction.
\qed

\begin{theo} \label{exists}
If there exists $\delta_Y$ such that for all $\delta< \delta_Y$
we have $\tilde{Y}^\delta=\tilde{Y}^{\delta_Y}$ then the universal
cover of $Y$ exists and is $\tilde{Y}^{\delta_Y}$.
\end{theo}

\Pf 
We need only show that given any cover, $\tilde{Y}'$, of $Y$, 
$\tilde{Y}^{\delta_Y}$ covers $\tilde{Y}'$.  By Theorem 12 in
[Sp, p81], $p:\tilde{Y}'\to Y$ is a covering projection only
if there is an open covering ${\cal U}$ of $Y$ and a point
$\tilde{x}\in \tilde{Y}'$ such that
$$
\pi({\cal U}, p(\tilde{x})) \subset p_{\#}\pi(\tilde{Y}', \tilde{x}).
$$
Since $Y$ is compact, there is a finite sub-cover, ${\cal U}'$, 
such that $\pi(Y,{\cal U}, p(\tilde{x}))=\pi(Y,{\cal U}', p(\tilde{x}))$.
For all $y\in Y$ there is a $r_y>0$ such that $B (y, r_y)\subset U$
and $U\in {\cal U}'$.  Let 
$$
{\cal V}=\{ B_y(r_y):  y \in Y \}.
$$
Now ${\cal V}$ also has a finite sub-cover ${\cal V}'$, and
taking 
$$
\delta_0=\min \{r_y: B_y(r_y)\in \cal{V'},  \delta_Y\},
$$
we have a $\delta_0$ open cover, ${\cal U}_{\delta_0}$, which
refines ${\cal U}$.  Thus by [Sp, p 81 st 8],
$$
\pi({\cal U}_{\delta_0}, p(\tilde{x})) 
\subset p_{\#}\pi(\tilde{Y}', \tilde{x}).
$$
In particular, $\tilde{Y}^{\delta_0}$ covers $\tilde{Y}'$.  However
$\delta_0 \le \delta_Y$, so  
$\tilde{Y}^{\delta_Y}=\tilde{Y}^{\delta_0}$ covers $\tilde{Y}'$.
\qed

We now prove that the revised fundamental group of the limit space
of a sequence of simply connected compact length spaces is trivial
[Theorem~\ref{simpconn}].
The example in Section 2 shows that this is as much as we can hope for.

\vspace{.2cm}
\noindent{\bf Proof of Theorem~\ref{simpconn}:}
Since $M_i$ are simply connected, for all $\delta>0$, we know
$G(M_i, \delta)$ is trivial.  Thus by Corollary~\ref{GG},
we know that for all $\delta>0$, $G(Y, \delta)$ is trivial as well.
Thus $\tilde{Y}^\delta=Y$ for  all $\delta$ and we satisfy the conditions
of Theorem~\ref{exists}.  So the universal cover exists and 
has a trivial group of deck transformations.
\qed

We now apply these results to study limits of Riemannian manifolds with 
a lower bound on Ricci curvature.

%--------------------------------------------------
\sect{Ricci Curvature}
%--------------------------------------------------------------------------

In this section we assume $Y$ is a Gromov-Hausdorff limit 
of compact manifolds $\{ (M_i^n, p_i) \}$ satisfying (\ref{cond1}), that
is $\Ric M_i^n \geq -(n-1)$ and $Diam (M_i) \leq D$.  Recall that,
by the Gromov Precompactness Theorem \cite[Thm 5.3]{Gr}, a subsequence
of any sequence of covering spaces of such manifolds converges in the
Gromov-Hausdorff sense to some limit space.  In particular, we know
the following.

\begin{lem}\label{precomp}
If $Y$ is the limit of $\{ (M_i^n, p_i) \}$ satisfying (\ref{cond1})
and $\delta>0$, then there exists a subsequence such that
$\{ (\tilde{M}_{i_j}^\delta, \tilde{p}_i ) \}$ converges in the
Gromov-Hausdorff sense to some limit space $Y^\delta$
with all the properties of Theorem~\ref{limdelcov}.
\end{lem}

Thus we can apply Cheeger-Colding's result \cite{ChCo} to show that 
Bishop-Gromov's volume comparison theorem holds on $Y^\delta$.
Recall \cite{ChCo}

\begin{theo}[Cheeger-Colding]
Given any sequence of pointed manifolds, $\{ (M_i^n, p_i) \}$, for 
which $\Ric M_i^n \geq -(n-1)$ holds, there is a subsequence, 
$\{ (M_j^n, p_j) \}$, convergent to some $(Y^m,y)$ in the pointed 
Gromov-Hausdorff sense, and there is a measure $\bar{V}_\infty$ 
on $Y$ satisfying Bishop-Gromov's volume comparison theorem, i.e. for 
$z \in Y^m$, $r_1 \leq r_2$, the following holds:
\be  \label{bishop}
\frac{\bar{V}_\infty(z, r_1)}{\bar{V}_\infty (z, r_2)} 
\geq \frac{V_{n,-1} (r_1)}{V_{n,-1} (r_2)},
\ee
where $V_{n,-1} (r_1)$ is the volume of a ball of radius $r_1$ in 
the simply connected space of dimension $n$ and curvature $\equiv -1$. 
\end{theo}

Therefore, we immediately have the following

\begin{coro} \label{vol}
Bishop-Gromov's volume comparison theorem (\ref{bishop}) holds on  
the Gromov Hausdorff limit, $Y^\delta$, of any converging subsequence of 
$\delta$-coverings of $M_i$.
\end{coro}

In \cite{ChCo}, Cheeger and Colding prove the following theorem
about the regularity of the limit spaces of spaces with Ricci curvature 
curvature bounded below. 

\begin{defn}{\em  
A {\em regular point}, $y$, in a limit space, $Y$, is a point such that
there exists $k$ such that every tangent cone at $y$ is isometric
to ${\Bbb R}^k$.}
\end{defn}

\begin{theo}[Cheeger-Colding] \label{regular} 
If $Y$ is the limit space of a sequence of $M_i^n$ with 
$\Ric M_i^n \geq -(n-1)$ then the set of
regular points has positive measure, $\bar{V}_\infty({\cal R})>0$.
In particular, the regular points are dense in $Y$.
\end{theo}

Although we cannot control the topology in general, near 
a regular point, we can control the $\delta$-covers.  
In fact we can control these $\delta$-covers above points with poles 
in all their tangent cones.  
Recall that a tangent cone $Y^\infty$, has a {\em pole} at
$y_\infty$ if  $\forall x \in Y^\infty$, there
is a length minimizing curve emanating from $y_\infty$ that
passes through $x$ and extends minimally to $\infty$.
Naturally this occurs at a regular point.

\begin{theo}  \label{regdel}
Let $Y$ be the Gromov-Hausdorff limit of a sequence of compact manifolds
such that
\be \label{cond3}
\Ric \ge -(n-1)K \textrm{ where } n\ge 3, K>0.
\ee
If
$y \in Y$ is a point such that there exists a tangent cone, 
$Y^\infty, y_\infty$,
that has a pole at $y_\infty$, then there exists $r_y>0$, such that for 
all $\delta>0$, $B(y,r_y)$ lifts isometrically to $Y^\delta$.
\end{theo}

The proof of this theorem uses an idea similar to one used in \cite{So}. 
Here, however, we have an arbitrary lower bound on Ricci curvature and are
concerned with eliminating small loops rather than large ones.  
In both cases we need to look at the Excess Theorem of Abresch 
and Gromoll from a new perspective \cite{AbGl}.  
This new perspective is required because the original excess
theorem has an inequality that includes the distance from a point
to a minimal geodesic.  Such a distance does not adapt well here
because a curve of minimal length in a limit space is not 
necessarily the limit of minimal geodesics.  Thus we begin with
the following lemma.

\begin{lem}  \label{LemSnk} \label{ptexcess}
Let $M^n$ be a complete Riemannian manifold satisfying
(\ref{cond3}).  There exists a constant
\be \label{EqnSnk}
S=S_{n, K}=\min \left\{  \frac{1}{8},
\frac{1}{4\cdot 3^{n}} \frac{1}{\cosh \sqrt{K}/4} \frac{n}{n-1}  
\left(\frac{n-2}{n-1} \right)^{n-1} 
\left( \frac{\sqrt{K}}{\sinh \sqrt{K}} \right)^{n-1} \right\}
\ee
such that if  $\gamma$ is a minimal geodesic of length $D\le 1$ 
and $x \in M^n$ satisfying
$$
d(x, \gamma(0))\ge (S_{n,K}+1/2)D \textrm{ and } 
d(x, \gamma(D))\ge (S_{n,K}+1/2)D. 
$$
then 
\be \label{abgl0}
d(x, \gamma(D/2)) \ge 3S_{n,K}D.
\ee
\end{lem}

Note that in the case where with nonnegative Ricci curvature, 
$D$ can have arbitrary length and $S$ has no hyperbolic terms \cite{So}.

\Pf
First we recall the Excess Theorem \cite[Prop. 2.3]{AbGl}.  If
$r_0=d(x, \gamma(0))$, $r_1=d(x, \gamma(D))$ and
$l=d(x, \gamma)<\min \{r_0, r_1\}$, then 
\be \label{abgl1}
e(x)=r_0+r_1-D \leq 2 \left( \frac{n-1}{n-2}\right) 
\left( \frac{1}{2} C_3 l^n \right)^{1/(n-1)},
\ee
where
\be
C_3 = \frac{n-1}{n} \left( \frac{\sinh \sqrt{K} l}{\sqrt{K} l} \right)^{n-1}
\sqrt{K} \left[ \coth \sqrt{K} (r_0 -l) + \coth \sqrt{K} (r_1 -l) \right].
\ee

Now suppose $l>\min \{r_0,r_1\}$, then 
$$
d(x, \gamma(D/2)) \ge l \ge (S+1/2)D \ge 3 SD.
$$
So we need only consider the case where  $l<\min \{r_0,r_1\}$
and we can apply (\ref{abgl1}).

Let us assume on the contrary that (\ref{abgl0}) does not hold.
Then $l < 3SD$,
$ r_0-l > (1/2+S)D-3SD > D/4$ and similarly $r_1-l > D/4$. Therefore 
\ban
e(x) & < 
& 2 \left( \frac{n-1}{n-2}\right) 
\left[ \frac{1}{2} \left(\frac{n-1}{n} \right) 
\left( \frac{\sinh (\sqrt{K} 3SD)}{\sqrt{K} 3SD} \right)^{n-1}  
2\sqrt{K} \coth (\sqrt{K} \frac{D}{4}) (3SD)^n \right]^{1/(n-1)} \\
& \leq &  2 \left( \frac{n-1}{n-2}\right) 
\left[ \frac{1}{2} \left(\frac{n-1}{n} \right) 
\left( \frac{\sinh \sqrt{K} }{\sqrt{K}} \right)^{n-1}  
2\sqrt{K} \frac{\cosh (\sqrt{K}/4)}{\sqrt{K} 
\frac{D}{4}} (3SD)^n \right]^{1/(n-1)} \\
& \leq & 2D  \left( \frac{n-1}{n-2}\right) 
\left[ 4 \left(\frac{n-1}{n} \right) 
\left( \frac{\sinh \sqrt{K} }{\sqrt{K}} \right)^{n-1} 
\cosh (\sqrt{K}/4) (3S)^n \right]^{1/(n-1)}.
\ean
On the other hand, $e(x)=r_0+r_1 -D \ge 2(S+1/2)D-D=2SD$, so
\[
S <  \frac{n-1}{n-2} \left[ 4 \left(\frac{n-1}{n} \right) 
\left( \frac{\sinh \sqrt{K} }{\sqrt{K}} \right)^{n-1} 
\cosh (\sqrt{K}/4) (3S)^n \right]^{1/(n-1)}.
\]
This gives 
\[
S > 4^{-1} 3^{-n} \frac{1}{\cosh \sqrt{K}/4} 
\frac{n}{n-1}  \left(\frac{n-2}{n-1} \right)^{n-1} 
\left( \frac{\sqrt{K}}{\sinh \sqrt{K} } \right)^{n-1},
\]
contradicting (\ref{EqnSnk}).
\qed 

We will now apply this lemma to prove Theorem~\ref{regdel}.  

\vspace{.2in}

\noindent{\bf Proof of Theorem~\ref{regdel}:}
Assume on the contrary that for all $r>0$ there is a $\delta_r>0$
such that the ball $B(y,r)$ does not lift isometrically to $Y^{\delta_r}$.   
Let $G^\delta$ denote the deck transformation group on $Y^\delta$. Thus,
there exist $r_i \to 0$, $\delta_i=\delta_{r_i}$,
and $g_i\in G^{\delta_i}$ such that 
$d_i=d_{{Y}^{\delta_i}}(\tilde{y}, g_i\tilde{y}) \in (0, 2r_i) \subset (0,1]$.
In fact, we can choose $g_i$ so that
\be \label{shortest}
d_{{Y}^{\delta_i}}(\tilde{y}, g_i\tilde{y})
\le d_{{Y}^{\delta_i}}(\tilde{y}, h\tilde{y}) \ \ 
\forall h\in  G^{\delta_i}.
\ee

Next we will find a length minimizing curve, $\tilde{C}_i$, running from 
$\tilde{y}$ to  $g_i\tilde{y}$ which has the property that it
passes through a particular point $\tilde{z}_i= \tilde{C}_i(d_i/2)$
which is the limit of halfway points of length minimizing curves in 
the sequence $\tilde{M_j}^{\delta_i}$.  We do this so that we can apply
Lemma~\ref{ptexcess} to $M_j$.  

To construct $\tilde{C}_i$, we first let $\tilde{y}_j$ and $\tilde{y}^i_j$ 
be points in $\tilde{M_j}^{\delta_i}$ which are close to $\tilde{y}$ and 
$g_i\tilde{y}$.  So $d_{\tilde{M_j}}(\tilde{y}_j, \tilde{y}^i_j)=d_{i,j}$
converges to $d_i$.  Let $\tilde{z}^i_j$ be midpoints of minimal geodesics
$\gamma^i_j$,
running from $\tilde{y}_j$ to $\tilde{y}^i_j$.  Taking a subsequence
of $j\to\infty$, there is a point $\tilde{z}_i\in Y^{\delta_i}$ which 
is halfway between $\tilde{y}$ to  $g_i\tilde{y}$.  Let $\tilde{C}_i$ 
be a length minimizing curve running from $\tilde{y}$ to $\tilde{z}_i$ 
and then to $g_i\tilde{y}$.  Finally let $C_i$ be the projection of 
$\tilde{C_i}$ to $Y$.

Now, imitating the proof of the Halfway Lemma of [So], 
and using (\ref{shortest}), we 
know $C_i\in Y$ is minimizing halfway around, $d_Y(C_i(0), C_i(d_i/2))=d_i/2$.

We choose a subsequence of these $i$ such that $(Y,y)$ rescaled
by $d_i$ converges to a tangent cone $(Y^\infty, y_\infty)$.  So
\begin{equation}  \label{tancone}
d_{GH}\left( B(y,10d_i)\subset Y, B(y_\infty, 10d_i)\right) < \epsilon_i d_i
\end{equation}
where $\epsilon_i$ converges to $0$.

Let $S$ be the constant from Lemma~\ref{LemSnk}.
Since $Y^\infty$ has a pole at $y_\infty$, we know there is a length
minimizing curve running from $y_\infty$ through any point in
$\partial B(y_\infty, d_i/2)$ to $\partial B(y_\infty, d_i/2+2S d_i)$. 
Thus by (\ref{tancone}), 
\be \label{wherex}
\forall \,\,\,x\in \partial B(y,d_i/2+2S d_i) \subset Y,
\ee
we have points
\be
x_\infty\in Ann_{y_\infty}(d_i/2+2S d_i-\epsilon_id_i,
d_i/2+2S d_i+\epsilon_id_i)
\ee
and
\be
y_i\in Ann_{y_\infty}(d_i/2-\epsilon_id_i,
d_i/2+\epsilon_id_i)
\ee
such that
\begin{eqnarray} \label{epsicont}
d_Y(x, C_i(d_i/2)) &<& d_{Y^\infty}(x_\infty, y_i) + \epsilon_i d_i\\
&\le & 2\epsilon_i d_i +2S d_i +\epsilon_i d_i.
\end{eqnarray}

Now we will imitate the Uniform Cut Lemma of [So], to show that
for all $x\in \partial B(y, d_i/2+2S d_i)$, we have
$l_i=d_Y(x, C_i(d_i/2)) \ge (3S)d_i.$
This will provide a 
contradiction for $\epsilon_i < S/2$ and we are done. 

First we lift our points $x$ and $y$ to the cover $\tilde{Y}^{\delta_i}$
as follows.  We lift $y$ to the point $\tilde{y}$ and we lift the closed
loop $C_i$ to the curve $\tilde{C}_i$ running from $\tilde{y}$
through $z_i=\tilde{C_i(d_i/2)}$ to $g_i\tilde{y}$.  Then if $\sigma$
is a length minimizing curve of length $l_i$ running from
$C_i(d_i/2)$ to $x$, we lift it to $\tilde{Y}^{\delta_i}$ so it
runs from $\tilde{z}_i$ to a new point, $\tilde{x}$.  
Note that by our choice of $x$ in (\ref{wherex}),
\be
d_{\tilde{Y}^{\delta_i}}(g_i\tilde{y}, \tilde{x})\ge
d_{{Y}}({y}, {x})=d_i/2+2S d_i
\ee
and so is $d_{\tilde{Y}^{\delta_i}}(\tilde{y}, \tilde{x})$.

By our choice of $\tilde{C}_i$ and $\tilde{z}_i$, we know there are 
corresponding
points in $\tilde{M}_j^{\delta_i}$.  That is there is a triangle
formed by $\tilde{y}_j$, $\tilde{y}^i_j$, 
with a minimal geodesic $\gamma^i_j$ running between them  and some point
$\tilde{x}_j$ such that 
\begin{eqnarray*}
d_{i,j}=d_{\tilde{M_j}^{\delta_i}}(\tilde{y}_j, \tilde{y}^i_j) & \to& d_i ,\\
d_{\tilde{M}_j^{\delta_i}}(\tilde{y}_j, \tilde{x}_j)&\to&
d_{\tilde{Y}^{\delta_i}}(\tilde{y}, \tilde{x})=(1/2+2S)d_i\\
d_{\tilde{M}_j^{\delta_i}}(\tilde{y}^i_j, \tilde{x}_j)&\to&
d_{\tilde{Y}^{\delta_i}}(g_i\tilde{y}, \tilde{x})=(1/2+2S)d_i.\\
l_{i,j}=d_{\tilde{M_j}^{\delta_i}}(\tilde{\gamma}^i_j(d_{i,j}/2), \tilde{x}_j)
&\to& d_{\tilde{Y}^{\delta_i}}(\tilde{z}_i, \tilde{x})=l_i.
\end{eqnarray*}
So for $j$ sufficiently large, we have
\be
d_{\tilde{M_j}^{\delta_i}}(\tilde{y}_j, \tilde{x}_j)\ge (1/2+S)d_{i,j}
\textrm{ and }
d_{\tilde{M_j}^{\delta_i}}(\tilde{y}^i_j, \tilde{x}_j)\ge (1/2+S)d_{i,j}
\ee
and can apply Lemma~\ref{ptexcess} to get
\be
l_{i,j} \ge 3Sd_{i,j}.
\ee

Taking $j$ to infinity, we get the limit of this bound in 
$\tilde{Y}^{\delta_i}$,  namely $l_i \ge 3S d_i$.
This contradicts (\ref{epsicont})
for $\epsilon_i < S/2$ and we are done. 
\qed

\begin{note}
In the non-collapsed case,  namely when the sequence of compact manifolds 
$M_i^n$ satisfy (\ref{cond2}),
every tangent cone of the limit space at every point is polar \cite{ChCo}. 
By Theorem~\ref{regdel}, for every $y \in Y$, there exists $r_y>0$, such 
that for all $\delta>0$,
$B(y,r_y)$ lifts isometrically to $Y^\delta$ for all $\delta >0$. Since 
$Y$ is compact, this implies that the $Y^\delta$ stabilize, i.e. there 
exists $\delta_Y$ depending on
$Y$ such that for all $\delta<\delta_Y$, we have 
$Y^\delta=Y^{\delta_Y}$. 
\end{note}

Using Theorem~\ref{regdel} and volume comparison we will prove that this 
is also true in the collapsed case. 

\begin{theo}  \label{samedelta}
There exists $\delta_Y$ depending on
$Y$ such that for all $\delta<\delta_Y$, we have 
$Y^\delta=Y^{\delta_Y}$ and $G^\delta=G^{\delta_Y}$. 
Therefore $\tilde{Y}^\delta$ are also same for all small $\delta >0$.
\end{theo}

\Pf Note that from Theorem~\ref{limdelcov}, 
if $Y^\delta$ do not stabilize for $\delta$ small, 
then neither do the $\delta$-covers $\tilde{Y}^\delta$.
So there exists a sequence of $\delta_i>0$ with 
$\delta_1 \leq D, \delta_i > 10 \delta_{i+1}$ 
such that all $\tilde{Y}^{\delta_i}$  and $G(Y,\delta_i)$ 
are distinct.   In particular there are elements of $G(Y, \delta_i)$
which are trivial in $G(Y, \delta_{i-1})$.  So there exist
$q_i\in Y$, such that the $B_{q_i}(\delta_{i-1})$ contains a noncontractible
loop, $C_i$, which lifts non-trivially in $Y^{\delta_i}$.  Since $C_i$
must lift to a union of balls $B_{g\tilde{q}_i}(\delta_{i-1})$ in
$\tilde{Y}^{\delta_i}$, there exists $g_i$ nontrivial in $G(Y, \delta_i)$
such that 
\be
d_{\tilde{Y}^\delta_i}(g_i\tilde{q}_i, \tilde{q}_i) < 2 \delta_{i-1}.
\ee 
So if $\alpha_i \subset Y$ is the projection of the minimal geodesic from
$g_i\tilde{q}_i$ to $ \tilde{q}_i$ it represents an element $g_i$
of ${\pi}_1(Y)$ which is mapped non-trivially into $G_{\delta_i}$ and
trivially into $G_{2\delta_{i-1}}$.

For any $j$, the limit cover $Y^{\delta_j}$ covers
$Y^{\delta_i}$ for $i=1... j-1$.  So $g_1...g_{j-1}$ are distinct 
nontrivial deck transforms of $Y^{\delta_j}$.  Furthermore, for any $q\in Y$, 
letting $\tilde{q_i}$ be the lift of $q_i$ closest to
 $\tilde{q}\in Y^{\delta_j}$,
we have,
\be
d_{Y^{\delta_j}}(\tilde{q}, g_i\tilde{q})
\le d_{Y^{\delta_j}}(\tilde{q}, \tilde{q_i})
+d_{Y^{\delta_j}}(\tilde{q_i}, g_i\tilde{q_i})
+d_{Y^{\delta_j}}(g_i\tilde{q_i}, g_i\tilde{q})
\le D + 2 L(\alpha_i) +D \le 4D.
\ee 
Therefore we have for any $j$, 
there are $j-1$  distinct elements in 
$G^{\delta_j}$ with $l(g_i, \delta_j) \leq 4D$.

On the other hand the total number of elements in $G^\delta$ of 
$\delta$-length $\leq 4D$ is uniformly bounded for all $\delta$ in 
terms of geometry and topology of $Y$. To show this let us look 
at the lift of a regular point $p \in Y$ in the cover $Y^\delta$.  We know
by Theorem~\ref{regdel}, there is a $\delta_0>0$ such that 
the ball of radius $\delta_0$ about $p$ is isometrically lifted 
to disjoint balls of radius $\delta_0$ in $Y^\delta$. 
Let $N$ be the number of distinct elements in $G^\delta$ of 
$\delta$-length $\leq 4D$. Note that $g B(\tilde{p}, \delta_0)$ 
is contained in $B(\tilde{p}, 4D+\delta_0 )$ for all $g \in G^\delta$ 
with $l(g, \delta) \leq 4D$.
Thus applying Corollary~\ref{vol} we have   
\be
N \leq  \frac{\bar{V}_\infty (\tilde{p}, 4D+\delta_0)}
{\bar{V}_\infty (\tilde{p}, \delta_0)} 
\leq  \frac{V_{n,-1} (4D+\delta_0)}{V_{n,-1} (\delta_0)}.
\ee

 This is a contradiction.
\qed

Note that $\delta_Y$ can not be uniformly bounded regardless of $Y$,
because we only bound the number of elements of a certain length.
Furthermore, the $\delta_0$ depended on the properties of the regular point
in $Y$.  Finally, even with a uniform lower bound on volume,
a uniform bound on $\delta_Y$ would contradict Otsu's examples.

This result has several nice consequences.  First, combining this with 
Theorem~\ref{exists}, we get the following.

\begin{theo} \label{ricciexists}
If $Y$ is the limit of a sequence of compact manifolds with 
uniformly bounded diameter and a uniform lower bound on Ricci
curvature, then the universal  covering space of $Y$ exists.
\end{theo}

\begin{note}
In fact $Y^{\delta_Y}$ is the universal cover.
It is unknown whether the universal cover is simply connected.
\end{note}

This theorem allows us to define the revised fundamental group
of $Y$, $\bar{\pi}_1(Y)$ [Definition~\ref{deck}].  We can now prove
Theorem~\ref{surjective} which was stated in the introduction.

\vspace{.2cm}
\noindent{\bf Proof of Theorem~\ref{surjective}:}
First, there is always a surjection 
\be
\pi_1(M_i, p) \to G(M_i,\delta)=\pi_1(M_i, p)/\pi_1(M_i, \delta, p).
\ee
Take $\delta=\delta_Y/2$.  Then by Theorem~\ref{theodelsurj},
there exists $N$ sufficiently large that there is a surjection
$\Phi_i: G(M_i,\delta) \to G(Y,\delta_Y)$.  Since
$\tilde{Y}^\delta$ is the universal cover of $Y$, 
$G(Y,\delta_Y)=\bar{\pi}_1(Y)$.
\qed

Theorems~\ref{ricciexists} and~\ref{limdelcov} imply the following.

\begin{coro} \label{subdelcov}
The universal cover of a limit space, $Y$, of compact manifolds, $M_i$, 
satisfying
(\ref{cond1}), is the limit of $\delta_Y/2$-covers of a subsequence
of the $M_i$.  
\end{coro}

\begin{coro}  \label{theoisom}
For all $\delta<\delta_Y$, there exists $N$ sufficiently large depending on 
$\delta$ and $\delta_Y$ such
that $G(M_i, \delta)=\bar{\pi}_1(Y)$ for all $i \ge N$.
\end{coro}

\Pf
Choose $N$ sufficiently large that 
$$
d_{GH}(M_i, Y) < 
\min \left\{\frac{\delta}{20}, \frac{\delta_Y-\delta}{20} \right\}
\qquad \forall i \ge N.
$$
Then  by Theorem~\ref{theodelsurj}, there exists a surjection
from $G(M_i, \delta)$ to $G(Y,\delta_Y)$ and there
is also a surjection from $G(Y,\delta/2)=G(Y,\delta_Y)$ to $G(M_i,\delta)$.
It is clear from the definition of these surjections in the proof
that they commute.  
\qed

In the non-collapsed case, i.e. the sequence of compact
manifolds $M_i^n$ satisfy (\ref{cond2}),
we can now prove Theorem~\ref{noncol} 
that the fundamental groups of $M_i$ are eventually
isomorphic up to a finite group.
 
\vspace{.2cm}
\noindent{\bf Proof of Theorem~\ref{noncol}:}
This essentially follows from
Corollary~\ref{theoisom} and
 Anderson's estimate on the order of subgroups generated by short 
loops \cite[Page 268]{An} (see also \cite[Page 256]{Pe2}).
Anderson's estimate states that there 
exist $L= L(n,v,D)$ and $N=N(n,v,D)$ such 
that if $M^n$ is a compact manifold satisfying (\ref{cond2}), then 
any subgroup of $\pi_1(M)$ that is generated by loops of length 
$\leq L$ must have order $\leq N$. Now choose $\delta$ such that 
$\delta \leq \min \{L/2, \delta_Y/2 \}$. Since the covering group 
of the $\delta$-cover of each $M_i$, $\pi_1(M_i, \delta)$, is 
generated by loops of length $\leq 2\delta \leq L$, the order 
of $\pi_1(M_i, \delta)$ is uniformly bounded by $N(n,v,D)$. On 
the other hand, by Corollary~\ref{theoisom}, the deck transformation 
group of the $\delta$-covering space of $M_i$, 
$G(M_i,\delta) = \pi_1(M_i)/\pi_1(M_i, \delta)$, is isomorphic to 
$\bar{\pi}_1(Y)$ for all $i \geq i_0$. 
Setting $F_i = \pi_1(M_i, \delta)$ finishes the proof.
\qed

Department of Mathematics and Computer Science, 

Lehman Colleger, City University of New York,

Bronx, NY 10468

sormani@g230.lehman.cuny.edu

Department of Mathematics, 

University of California, 

Santa Barbara, CA 93106 

wei@math.ucsb.edu

\end{document}